
\input amstex
\documentstyle{amsppt}
\magnification\magstep1
\NoBlackBoxes

\def\R{{\Bbb R}}
\def\Z{{\Bbb Z}}

\def\Rp#1{\R\roman P^{#1}}
\def\Ker{\mathop{\roman{Ker}}\nolimits}

\pagewidth{28pc}
\pageheight{43pc}
\topskip=\normalbaselineskip


\let\tm\proclaim
\let\endtm\endproclaim


\let\pf=\demo
\let\endpf=\enddemo


\let\rk=\remark
\let\endrk=\endremark

\let\de=\definition
\let\endde=\enddefinition


\begingroup
\let\head\relax \let\subhead\relax 

\outer\gdef\ah#1\endah{\head#1\endhead}

\outer\gdef\bh#1\endbh{\subhead#1\endsubhead}

\endgroup 

\def\a{\alpha}
\def\b{\beta}
\def\om{\omega}
\def\F{\Phi}
\def\Pin{\operatorname{Pin}}
\def\Spin{\operatorname{Spin}}
\def\O#1{\widetilde O_{#1}}
\def\marginalnote#1{}
\def\<#1,#2>{\left<#1,#2\right>}
\def\tmod#1{$\operatorname{mod}#1$}

\topmatter
\title
$\Pin$-structures on surfaces and quadratic forms
\endtitle
\author
A.~Degtyarev, S.~Finashin
\endauthor
\address
UC Riverside, Riverside, CA 92521
\endaddress
\email
degt\,\@\,ucrmath.ucr.edu
\endemail
\address
Middle East Technical University, Ankara, 06531, Turkey
\endaddress
\email serge\,\@\,rorqual.cc.metu.edu.tr
\endemail
\abstract
A correspondence between different $\Pin$-type structures on  a compact
surface and quadratic (linear) forms on its homology is constructed.
Addition of structures is defined and expressed in terms of these
quadratic forms.
\endabstract
\endtopmatter

\document
In this paper we try to clarify the relation between $\Pin$-type
structures on a compact surface and quadratic forms on its homology.
This relation is well-known and useful in the case of $\Spin$- and
$\Pin^-$-structures (\cite{J}, \cite{KT}\marginalnote{cites J, KT}).
It makes it much easier to understand the nature of some fundamental
in low-dimensional topology objects such as $\Z/2$-Seifert form on a
surface in an~oriented 3-manifold and Rokhlin form on
a~characteristic surface in an oriented 4-manifold (see, e.g.,
\cite{F},~\cite{DFM}).

\rk{Remark}
A slightly more general approach of \cite{DFM} also explains these
forms in the case of a non-oriented ambient manifold, as well as the
newly found Benedetti-Marin form~\cite{BM}.
\endrk

We will show that a similar correspondence between quadratic forms
and structures can also be defined for other $\Pin$-type structures.
In this short paper we restrict ourselves to a geometrical description
in the simplest case, when
the structural group is a $\Z/2$-extension of the orthogonal
group~$O_n$. Up to isomorphism, there are four such extensions
corresponding to the four elements of~$H^2(BO_n;\Z/2)$, each element
being the obstruction to existence of a structure in
an~$O_n$-bundle $P\to X$. When non-empty, the set of all the
structures of a given type forms an affine space over~$H^1(X;\Z/2)$.
If the obstruction is the trivial element of $H^2(BO_n;\Z/2)$,
this set obviously
coincides with $H^1(X;\Z/2)$. The classes $w_2$ and $w_2+w_1^2$
characterize $\Pin^+$- and $\Pin^-$-structures respectively.
Finally, the structures corresponding to the remaining class $w_1^2$
are not (to the best of our knoledge) mentioned in literature
and did not receive any special name.
We will call them $\O n$-structures, $\O n$ standing for the
nontrivial semi-direct product $\Z/4\ltimes SO_n$
(topologically, $\O n\to O_n$ is a trivial double covering).

\rk{Remark}
Actually, $\O n$-structures do appear in literature implicitly,
e.g., as framings in the complexification of a vector bundle~\cite{A},
or as linear forms on a real algebraic variety~\cite{N}. Besides,
they complete the descending table in~\cite{KT, Corol\.~2.15}:
a~$\Pin^-$-structure on a manifold $M$ descends to an
$\O n$-structure on a codimension one submanifold whose normal
bundle is isomorphic to the determinant of the tangent bundle
of~$M$.
\endrk

Thus, we show that for each of these four classes there is a
one-to-one correspondence between the set of structures on a compact
surface and the set of specific quadratic (or, in special cases,
linear) forms on the 1-homology of the surface. (This correspondence
is known in the $\Pin^-$-case (see, e.g.,~\cite{KT}) and is
obvious in the case of the trivial extension. The two others are
defined in \S2.) Another subject of the paper, which has never (as
far as we know) been mentioned explicitly (see, though, a slightly
different approach in~\cite{D}), is addition of~structures. This
operation, defined in \S3, naturally extends the canonical affine
action of~$H^1(X;\Z/2)$ and covers addition of the characteristic
classes. (By the way, this gives one more reason for considering
$\O n$-structures: they are sums of $\Pin^-$ and $\Pin^+$-structures.)
We give an interpretation of this
operation in terms of quadratic forms.


\ah
\S2. Quadratic forms
\endah

\bh
1.\enspace $\Pin^-$-structures---quadratic forms $H_1(F;\Z/2)\to\Z/4$
\endbh
We start with reminding the standard construction of the quadratic
form $q$ corresponding to a $\Pin^-$-structure on a compact
surface~$F$.  Pick an integral class $\a\in H_1(F;\Z)$ and realize it
by an immersed collection of oriented circles $S\to F$. Let $n(S)$ and
$i(S)$ be the numbers of components and self-intersection points
of~$S$ respectively. The tangent vector field to $S$ defines a
$\Pin_1^-$-reduction of the restriction to $S$ of the given
$\Pin_2^-$-bundle on~$F$. Since $\Pin_1^-\cong\Z/4$ is a discrete
abelian group, one can consider the {\it total holonomy\/}
$h(S)\in\Pin_1^-$ of this bundle along~$S$ (which, by definition, is
the sum in $\Pin_1^-$ of the holonomies along the components of~$S$).
We let $q(\a)=h(S)+2\bigl(n(S)+i(S)\bigr)\pmod4$. Now standard
arguments apply to show that $q(\a)$ does not depend on $S$ and
satisfies the identity $q(\a+\b)=q(\a)+q(\b)+2\<\a,\b>$, where
$\<\cdot,\cdot>$ is the intersection form on~$F$ and
$2\:\Z/2\to\Z/4$ is the unique inclusion. ($q(\a)$ obviously does not
change during a regular homotopy of~$S$, and elementary
transformations like Reidemeister move I and smoothing a
self-intersection point can easily be controlled; see, e.g.,
\cite{KT}.) Since the \tmod2 reduction of~$q$ coincides
with $w_1\:H_1(F;\Z/2)\to\Z/2$, the above formula implies, in
particular, that $q$ factors through $\Z/2$-homology of~$F$.

\bh
2.\enspace $\O2$-structures---linear forms $H_1(F;\Z/4)\to\Z/4$
\endbh
Since the corresponding 1-dimensional group $\O1$ is also $\Z/4$,
this case is similar to the previous one. The only difference is that
one should not adjust holonomy by the numbers of components and
self-intersection points, i.e., $q(\a)=h(S)$. The result is a linear
form $q$ which factors through $\Z/4$-homology,
$q\:H_1(F;\Z/4)\to\Z/4$, and whose restriction \tmod2 coincides with
$w_1\:H_1(F;\Z/2)\to\Z/2$.

\bh
3.\enspace $\Pin^+$-structures---quadratic forms $H_1(F;\Z/4)\to\Z/2$
\endbh
The construction goes similar to the case of $\Pin^-$-structures.
Now $\Pin_1^+\cong O_1\times\Z/2$, the projection of~$h(S)$ to the
first factor $O_1$ being just the value of $w_1$ on $\a$. In order to
drop this standard component, we consider the projection $p_2h(S)$
to the second factor~$\Z/2$; then we let
$q(\a)=p_2h(S)+n(S)+i(S)\pmod2$. This form satisfies the identity
$q(\a+\b)=q(\a)+q(\b)+\<\a,\b>$, which, in particular, implies
that it factors through $H_1(F;\Z/4)$.

\bh
4.\enspace Trivial structures---linear forms $H_1(F;\Z/2)\to\Z/2$
\endbh
This case, when structures just {\it are\/} cohomology classes, admits
a description similar to the previous three: the total holonomy $h(S)$
is an element of $O_1\times\Z/2$, and we let $q(\a)=p_2h(S)$.

\medpagebreak

We can now uniformize all the four cases and consider quadratic
(linear) forms $H_1(F;\Z/4)\to\Z/4$. (In the case of $\Pin^+$- and
trivial structures $\Z/2$ is embedded in $\Z/4$ via multiplication
by~2.) This gives the following result:

\tm{Theorem A}
Given a compact surface $F$, there is a canonical affine one-to-one
correspondence between structures on~$F$ with the characteristic class
$a w_2+b w_1^2$ \rom(for some fixed $a,b\in\Z/2$\rom) and functions
$q\:H_1(F;\Z/4)\to\Z/4$ satisfying the following conditions:
\roster
\item $q(\a+\b)=q(\a)+q(\b)+2a\<\a,\b>$;
\item $q(\a)=bw_1(\a)\pmod2$.
\endroster
\endtm

\demo{Proof}
The only thing that needs proof is the fact that the constructed map
$\{\text{structures}\}\to\{\text{forms}\}$ is one-to-one. Since
both the sets are affine over $H^1(F;\Z/2)$, it suffices to show that
existence of forms implies existence of structures. This is obvious
for $\Pin^-$- and trivial structures, or if the surface is not closed
(since structures always exist in such cases). For $\Pin^+$-
and $\O2$-structures on closed surfaces one can easily see that
desired forms exist if and only if elements of order~2 in
$H_1(F;\Z/4)$ annihilate $w_1$ (or, equivalently, have trivial
self-intersection).  This is the case when the surface is the
connected sum of an even number of $\Rp2$'s, i.e., exactly when
$w_2=w_1^2=0$.  \qed \enddemo

\tm{Corollary \rm(classification of $\Pin^+$-structures up to
isomorphism)}
Two $\Pin^+$-structures on a closed surface $F$ are isomorphic\/
\rom(i.e., can be transformed into each other by a diffeomorphism
of the surface\/\rom) if and only if the values of the corresponding
quadratic forms on the\/ \rom(unique\/\rom) \rom2-torsion element
of $H_1(F;\Z)$ coincide. In particular, two structures are isomorphic
if and only if they are cobordant.
\endtm

\demo{Proof}
The mentioned value is the only algebraic invariant of forms (an easy
exercise), and, as usual in 2-dimensional topology, one can find
an~automorphism of the lattice $H_1(F;\Z)$ which is accompanied by
a~diffeomorphism of~$F$.
\qed
\enddemo

\rk{Remark}
Note that we have to consider integral homology here, since otherwise
one cannot distinguish between different 2-torsion elements, and
the algebraic invariant disappears.
\endrk

\ah
\S3. Addition of structures
\endah

Given two $\Z/2$-extensions $G_1\to O_n$, $G_2\to O_n$, one can
define their {\it sum\/} $G_1\vee G_2$ to be the quotient
$G_1\times_{O_n}G_2\big/\operatorname{Diag}(\Z/2)$, where
$\text{Diag}$ is the canonical diagonal map
$$
\text{Diag}\:\Z/2=\Ker[G_i\to O_n]@>>>\Z/2\oplus\Z/2=
\Ker[G_1\times_{O_n}G_2\to O_n].
$$
(In fact, this is one of the standard algebraic approaches to
definition of the group structure on the set of isomorphism classes
of $\Z/2$-extensions of $O_n$, which is isomorphic to
$H^2(BO_n;\Z/2)$.) To apply this procedure to structures, one should
fix first some representatives $G(\om)$ of the isomorphism classes of
extensions, one for each characteristic class
$\om\in H^2(BO_n;\Z/2)$, and some maps
$G(\om_1)\vee G(\om_2)\to G(\om_1+\om_2)$. To do that uniformly in
all dimensions, it suffices to just pick some isomorphisms
$\Pin^-_1=\O1=\Z/4$, \ $\Pin^+_1=O_1\times\Z/2$, and
$\Z/4\vee\Z/4=O_1\times\Z/2$ (see \cite{DFM}). Such maps can
certainly be chosen and fixed once and forever. Then one can give the
following definition:

\de{Definition}
Let $P\to X$ be an $O_n$-bundle. Then, given two structures
$\F_1\to P$, $\F_2\to P$ with characteristic classes
$\om_1,\om_2\in H^2(BO_n;\Z/2)$ respectively, we define their {\it
sum\/} $\F_1\vee\F_2\to P$ to be the $(\om_1+\om_2)$-structure
associated with the fibered product $\F_1\times_P\F_2\to P$ via the
composed map
$$
G(\om_1)\times_{O_n}G(\om_2)@>>>G(\om_1)\vee G(\om_2)
  @>\cong>> G(\om_1+\om_2).
$$
\endde

Theorem B below is proved in \cite{DFM}.

\tm{Theorem B}
$\vee$ is a group operation on the set of all structures on a given
$O_n$-bundle $P\to X$, which extends the canonical affine action of
$H^1(X;\Z/2)$ on this set\/ \rom(i.e., $\vee$-sum with an
$(O_n\times\Z/2)$-structure coincides with the affine shift by the
corresponding cohomology class\/\rom).
\endtm

\tm{Theorem C}
$\vee$-sum of structures on a compact surface corresponds to the
following pointwise operation of quadratic forms:
$(q_1,q_2)\mapsto q_1+q_2+2q_1q_2$.
\endtm

\pf{Proof}
Due to the uniform construction of \S2 it suffices to consider
structures on the normal bundle to a circle, when the statement is
obvious.
\qed \endpf

The introduced $\vee$-sum operation admits an interpretation in
terms of $\Spin$-structures. Given an $O_n$-bundle $\xi\:P\to X$, let
us denote by $\Spin(\xi)$, $\Pin^\pm(\xi)$, etc\. the set of all the
$\Spin$-, $\Pin^\pm$-, etc\. structures on $\xi$ respectively. Then,
according to~\cite{KT}, there are natural isomorphisms
$\Pin^-(\xi)=\Spin(\xi\oplus\det\xi)$ and
$\Pin^+(\xi)=\Spin(\xi\oplus 3\det\xi)$. Similar arguments show that,
besides, there are isomorphisms
$\O n(\xi)=\Spin(2\xi)=\Spin(2\det\xi)$. Consider the Whitney sum of
the above three bundles:
$$
(\xi\oplus\det\xi)\oplus(\xi\oplus 3\det\xi)\oplus(2\xi)=
 4(\xi\oplus\det\xi).
$$
This bundle has a canonical $\Spin$-structure (the quaternion
$\Spin$-structure, which is defined on $4\eta$ for any bundle $\eta$,
see~\cite{DFM}). Hence, $\Spin$-structures on any two of the three
summands define a $\Spin$-structure on the third one, and one can
easily see that the obtained maps
$\Pin^-(\xi)\times\Pin^+(\xi)\to\O n(\xi)$, etc\. coincide with
the $\vee$-sum. This gives an alternative description of this
operation
in the most interesting cases which are not reduced to the affine
action of $H^1(X;\Z/2)$.

\enddocument